\documentclass[12pt]{article}
\usepackage{amsfonts,amssymb}
\usepackage{amsmath}
\usepackage{amsthm}
\usepackage[matrix,arrow,curve]{xy}
\newtheorem{theorem}{Theorem}
\newtheorem{lemma}{Lemma}
\newtheorem{prop}{Proposition}
\newtheorem{corollary}{Corollary}
\newtheorem{definition}{Definition}

\title{Tangent Hyperplanes to Subriemannian Balls}
\author{A.~A.~Agrachev\thanks{SISSA, Trieste \& Steklov Math. Inst.,
Moscow. The author has been supported by the grant of the Russian
Federation for the state support of research, Agreement No 14.b25.31.0029.}}
\date{}
\begin{document}
\maketitle
\begin{abstract}
We examine the existence of tangent hyperplanes to subriemannian balls. Strictly abnormal shortest paths are allowed.
\end{abstract}

\section{Introduction}

We still know very little on the balls and spheres in general (in fact, even in generic) subriemannian spaces.
In particular, we do not know if small balls are always contractible and small spheres are $C^0$-manifolds. The
difficulties are due to the presence of so called abnormal or singular geodesics connecting the center of the ball with the sphere. The distance is never Lipschitz in the points of these geodesics where it has rather mysterious singularities. Nevertheless, there is a good news at least under a generic assumption that all nonconstant geodesics have corank 1 (see Definition~3).
Let $\gamma(t),\ 0\le t\le r,$ be the parameterized by the length shortest path between $\gamma(0)$ and $\gamma(r)$. Theorem~1 of this paper states that centered at $\gamma(0)$ radius $t$ ball has a tangent hyperplane
at $\gamma(t)$, for any $t\in(0,r)$; see Definition~2 for the precise definition of the tangent hyperplane.

We do not distinguish normal and abnormal geodesics in the statement and in the proof of the theorem but actual geometry is very different. The statement is well-known for strictly normal geodesics; moreover, the radius $t$ sphere is smooth
in a neighborhood of $\gamma(t)$ in this case and is transversal to $\gamma$. The situation is much more delicate for abnormal geodesics: they are tangent to the spheres and this is why the distance is not Lipschitz.
A similar result with additional assumptions was proved in \cite{T}. Important extra assumptions were Whitney stratifiability of the sphere and uniqueness of the shortest path.

Theorem 2 of this paper is a counterpoint to Thorem~1: we show that any ball contains a lot of non-extendable
shortes paths such that the passing through their endpoints spheres do not admit tangent hyperplanes in these points.

\section{Subriemannian structures}

We start from basic definitions: see \cite{ABB} for more details and examples. Let $M$ be a smooth manifold. A subriemannian structure on $M$ is a homomorphism $f:U\to TM$, where $U\stackrel{\pi}{\longrightarrow}M$ is a smooth vector bundle equipped with an Euclidean structure. It is also assumed that $f$ is {\it bracket generating}.

Let us explain the last property that is also called the {\it H\"ormander condition}.
Let $\overline{U}$ be the $C^\infty(M)$-module of smooth sections of the bundle $U$, then $f(\overline{U})$ id a submodule of $\mathrm{Vec}M=\overline{TM}$. We say that a submodule of $\mathrm{Vec}M$ is involutive if it is a Lie subalgebra of $\mathrm{Vec}M$.
\begin{definition} A homomorphism $f:U\to TM$ is bracket generating if $f(\overline{U})$ is not contained in a proper involutive submodule of $\mathrm{Vec}M$.
\end{definition}

A Lipschitzian curve $\gamma:[0,t_1]\to  M$ is called horizontal if there exists a measurable bounded map $u:[0,t_1]\to U$ such that $\dot\gamma(t)=f(u(t))$ for a.\,e. $t\in [0,t_1]$. The subriemannian length of the horizontal curve is defined as follows:
$$
length(\gamma)\doteq\inf\left\{\int\limits_0^{t_1}|u(t)|_{\gamma(t)}dt : \dot\gamma(t)=f(u(t))\ \mathrm{for\ a.\,e.}\ t\in[0,t_1]\right\},
$$
where $|\cdot|_q$ is the Euclidean norm on the fiber $U_q=\pi^{-1}(q)$.

It is easy to see that any Lipschitzian reparameterization of a horizontal curve is a horizontal curve of the same length. Moreover, any horizontal curve is a reparameterization of a horizontal curve $\gamma$ parameterized by the length such that $length (\gamma\bigr|_{[0,t]})=t,\ \forall t\in[0,t_1]$.
The subriemannian or Carnot--Caratheodory distance between points $q_0,q_1\in M$ is defined as follows:
$$
\delta(q_0,q_1)=\inf\{length(\gamma): \gamma\ \mathrm{is\ a\ horizontal\ curve},\ \gamma(0)=q_0,\ \gamma(t_1)=q_1\}.
$$

Classical Chow--Rashevskij theorem implies that $\delta(\cdot,\cdot)$ is a well-defined metric on $M$ and this metric induces standard topology on $M$. In particular, small subriemannian balls are compact. Moreover, the subriemannian metric space is complete if and only if all balls are compact.

Given $q_0\in M,\ r>0$, we denote the centered at $q_0$ radius $r$ ball and sphere by $B_{q_0}(r)$ and $S_{q_0}(r)$, so that
$$
 B_{q_0}(r)=\{q\in M:\delta(q_0,q)\le r\},\quad S_{q_0}(r)=\{q\in M:\delta(q_0,q)=r\}.
$$
If $B_{q_0}(r)$ is compact then any point of $B_{q_0}(r)$ is connected with $q_0$ by a shortest horizontal path.
In other words, for any $q\in B_{q_0}(r)$ there exists a horizontal curve $\gamma$ such that
$\gamma(0)=q_0,\ \gamma(1)=q$ and $\delta(q_0,\gamma(1))=length(\gamma)$.

Any segment of a shortest path is, of course, a shortest path between its endpoints. In particular, if $\gamma$ is parameterized by the length and $\gamma(r)\in S_{\gamma(0)}(r)$, then $\gamma(t)\in S_{\gamma(0)}(t),\ \forall t\in[0,r]$.

Let $q\in M$ and $E_0\subset T_qM$ be a co-oriented hyperplane (vector subspace of codimension 1); then $T_qM\setminus E_0=E_+\cup E_-$, where $E_+$ and $E_-$ are positive and negative open half-spaces of $T_qM$.
\begin{definition} We say that a co-oriented hyperplane $E_0\subset T_qM$ is tangent to the ball $B_{q_0}(r)$ if for any smooth curve $\varphi:(-1,1)\to M$ such that $\varphi(0)=q$ and $\dot\varphi(0)\in E_+$ there exists $\varepsilon>0$ such that $\varphi(t)\notin B_{q_0}(r),\ \forall t\in(0,\varepsilon)$ and
$\varphi(t)\in B_{q_0}(r),\ \forall t\in(-\varepsilon, 0)$.
\end{definition}
Obviously, the ball may have at most one tangent hyperplane at the given point.

\section{The endpoint map}

Subriemannian structures $f_i:U_i\to M,\ i=1,2$ are called equivalent if there exists a commutative diagram
$$
\xymatrix{&U_1\ar[dr]_{f_1}&\\
V\ar[ur]_{g_1}\ar[dr]^{g_2}&&TM\\
&U_2\ar[ur]^{f_2}&}
$$
where $V$ is an Euclidean vector bundle over $M$, $g_1,g_2$ are epimorphisms and
$$
|u_i|=\min\{|v|: v\in V,\ g_i(v)=u_i\},\quad \forall u_i\in U_i,\ i=1,2.
$$
Equivalent structures have the same space of horizontal curves and the length of these curves is the same.

We say that $f:U\to TM$ is a free subriemannian structure if $U=M\times\mathbb R^k$, where $\mathbb R^k$ is equipped with the standard Euclidean inner product. Any subriemannian structure is equivalent to a free one.

Let $f:M\times\mathbb R^k\to TM$ be a free subriemannian structure, $e_1,\ldots,e_k$ be the standard orthonormal basis of $\mathbb R^k$ and $X_i(q)=f(q,e_i)$. A Lipschitzian curve $\gamma:[0,t_1]\to M$ is horizontal for this subriemannian structure if and only if there exist measurable bounded functions $u_i:[0,t_1]\to\mathbb R$ such that $\dot\gamma(t)=\sum\limits_{i=1}^ku_i(t)X_i(\gamma(t))$, for a.\,e. $t\in[0,t_1]$; moreover,
$length(\gamma)=\int\limits_0^{t_1}(\sum\limits_{i=1}^ku_i^2(t))^{\frac 12}dt$.

Now fix $q_0\in M$, a locally integrable vector-function 
$$
t\mapsto u(t)=(u_1(t),\ldots,u_k(t)),\quad t\ge 0, 
$$
and consider the Cauchy problem
$$
\dot q=\sum\limits_{i=1}^ku_i(t)X_i(q),\quad q(0)=q_0. \eqno (1)
$$
The problem has a unique Lipschitzian solution defined on an interval $[0,T)$. Let $t_1\in (0,T)$; then the solution to the Cauchy problem 
$$
\dot q=\sum\limits_{i=1}^kv_i(t)X_i(q),\quad q(0)=q_0,
$$ 
is defined on an interval that contains the segment $[0,t_1]$ for any $v$ sufficiently closed to $U$ in the $L_1$-norm.

In what follows, we prefer a more convenient $L_2$-norm. Let $\Omega^{t_1}_{q_0}$ be the set of all $u\in L_2^k[0,t_1]$ such that the solution of (1) is defined on the segment $[0,t_1]$. Then $\Omega_{t_1}$ is an open subset of $L_2^k[0,t_1]$.

The {\it endpoint map} $F^t_{q_0}:\Omega^t_{q_0}\to M$ sends $u\in\Omega_t$ to $q(t)$, where
$$
\dot q(\tau)=\sum\limits_{i=1}^ku_i(\tau)X_i(q(\tau)),\quad 0\le\tau\le t,\ q(0)=q_0. \eqno (2)
$$
Standard well-posedness results for the Cauchy problem imply that $F^t_{q_0}$ is a smooth map.

We have:
$$
\delta(q_0,q_1)=\inf\left\{\int\limits_0^t|u(\tau)|\,d\tau: u\in\Omega_t,\ F_t(u)=q_1\right\},\quad \forall t>0.
$$
Moreover,
$$
\delta^2(q_0,q_1)=\inf\left\{t\int\limits_0^t|u(\tau)|^2d\tau: u\in\Omega_t,\ F_t(u)=q_1\right\}, \eqno (3)
$$
as easily follows from the Cauchy--Schwartz inequality.

Now assume that $u\in\Omega^t_{q_0},\ \int\limits_0^t|u(\tau)|^2d\tau=t$ and $F_t(u)\in S_{q_0}(t)$. In other words, we assume that the solution of (2) is a shortest path connecting $q_0$ with $F_t(u)$. The relation (3) implies that $u$ is a critical point of the map $F^t_{q_0}$ restricted to the intersection of $\Omega_t$ with the $L_2$-sphere
$\{v:\int\limits_0^t|v(\tau)|^2d\tau=t\}$. The tangent space to the $L_2$-sphere at $u$ is equal to the
$u^\perp=\{v:\int\limits_0^t\langle v(\tau),u(\tau)\rangle\,d\tau=0\}$.

\begin{definition} A non-constant horizontal curve $\gamma$,
$$
\dot\gamma(\tau)=\sum\limits_{i=1}^ku_i(\tau)X_i(\gamma(\tau)),\quad 0\le\tau\le t,\  \gamma(0)=q_0, \eqno (4)
$$
is a subriemannian geodesic if $D_uF_t(u^\perp)\ne T_{\gamma(t)}M$. The codimension of the subspace $D_uF_t(u^\perp)\subset T_{\gamma(t)}M$ is called the corank of the geodesic $\gamma$.
\end{definition}

We explained that any shortest path is a geodesic. It is not hard to show that the corank depends only on $\gamma$ and not on the particular choice of $u(\cdot)$ that satisfies (4).

\medskip\noindent{\bf Example}. Let $M=\mathbb R^3=\{(x_1,x_2,x_3): x_i\in\mathbb R,\ i=1,2,3\},\ k=2,$
$$
X_1(x_1,x_2,x_3)=(1,0,0),\ X_2(x_1,x_2,x_3)=(0,1,x_1^n),\ q_0=(0,0,0,); 
$$
then
$$
F^1_{q_0}(u_1,u_2)=\left(\int\limits_0^1u_1(t)\,dt,\,\int\limits_0^1u_2(t)\,dt,\,
\int\limits_0^1u_2(t)\Bigl(\int\limits_0^tu_1(\tau)\,d\tau\Bigr)^ndt\right).
$$
The curve $t\mapsto F^t_{q_0}(0,1)$ is a geodesic of this structure. The corank of this geodesic equals 1 if $n=1$ and equals 2 if $n>1$.

\medskip
Let $\mathcal F(U,M)$ be the space of subriemannian structures $f:U\to M$ endowed with the $C^\infty$ Whitney topology and
$$
\mathcal F^1(U,M)=\{f\in\mathcal F(U,M):\ \mathrm{any\ geodesic\ of}\ f\ \mathrm{has\ corank\ 1}\}.
$$
It follows from the results of \cite{CJT} that $\mathcal F^1(U,M)$ contains an open dense subset of $\mathcal F(U,M)$. Moreover, $\mathcal F(U,M)\setminus \mathcal F^1(U,M)$ has codimension $\infty$ in $\mathcal F(U,M)$;
see \cite[Sec.\,5]{AG} for the precise definition of subsets of codimension $\infty$ and for the techniques that allow to easily prove this property.

\section{Main results}

Everywhere in this section we assume that $f\in\mathcal F^1(U,M)$ and the ball $B_{q_0}(r)$ is compact. We set:
$$
\Delta^1_q=\{f\bar u)(q):\bar u\in\bar U\},\quad \Delta^2_q=span\{[f(\bar u),f(\bar v)](q):\bar u,\bar v\in\bar U\}
$$

\begin{theorem} Let $q\in S_{q_0}(r)$ and $\gamma:[0,r]\to M$ a parameterized by the length shortest horizontal path connecting $q_0$ with $q$ s.\,t. $\gamma(0)=q_0,\ \gamma(r)=q$. Then $\forall\,t\in(0,r)$ there exists a tangent hyperplane to $S_{q_0}(t)$ at the point $\gamma(t)$.
\end{theorem}

\noindent{\bf Remark.} The corank 1 condition is essential. Indeed, the simplest example of a corank $>1$ shortest path is one from the previous section with $n=2$. The balls for this example are described in \cite{ABC};
they have acute angles at the points of the corank 2 geodesic and do not admit tangent planes in these points.

\begin{theorem} Assume that $\Delta^1_{q_0}\ne\Delta^2_{q_0}$. Then for any $t$ from an everywhere dense subset of $(0,r)$ there exists $q_t\in S_{q_0}(t)$ such that $S_{q_0}(t)$ does not admit a tangent hyperplane at $q_t$.
\end{theorem}

\noindent{\bf Proof of Theorem~1.} We have: $\gamma(\tau)=F^\tau_{q_0}(u)$, where
$|u(\tau)|\equiv 1,\ 0\le\tau\le r$. We fix $t\in(0,r)$ and denote $E^t=D_uF^t_{q_0}(u^\perp)$. We are going to prove that $E^t$ is the tangent hyperplane to $B_{q_0}(t)$ at $\gamma(t)$.

\begin{lemma} Let $\hat\gamma:[0,t]\to M$ be a horizontal path of the length $t$ with the same endpoints as $\gamma$:
$$
\hat\gamma(\tau)=F^\tau_{q_0}(\hat u),\ 0\le\tau\le t,\quad \hat\gamma(0)=q_0,\ \hat\gamma(t)=\gamma(t).
$$
 Then $D_vF^t_{q_0}(v^\perp)=E^t$.
 \end{lemma}

\noindent{\bf Proof of Lemma~1.} Let $P_t^r:M\to M$ be the diffeomorphism defined by the relation
$P_t^r:q(t)\mapsto q(r)$, where $q(\tau),\ t\le\tau\le r$, is a solution of the differential equation
$\dot q=\sum\limits_{i=1}^ku_i(\tau)X_i(q)$. In particular, $P_t^r(\gamma(t))=\gamma(r)$. The curve
$$
\beta(\tau)=\left\{\begin{array}{c}\hat\gamma(\tau),\ 0\le\tau\le t\\ \gamma(\tau),\ t\le\tau\le r
\end{array}\right.
$$
is a shortest path. Let $\hat\gamma(\tau)=F^\tau_{q_0}(\hat u),\ 0\le\tau\le t$, then
$\beta(\tau)=F^\tau_{q_0}(v),\ 0\le\tau\le r,$ where
$$
v(\tau)=\left\{\begin{array}{c}\hat u(\tau),\ 0\le\tau\le t\\ u(\tau),\ t\le\tau\le r
\end{array}\right.
$$
Hence $D_vF^r_{q_0}(v^\perp)$ is a codimension 1 subspace of $T_qM$. I claim that
$$
D_vF^r_{q_0}(v^\perp)=D_uF^r_{q_0}(u^\perp)=D_{u_t}F_{\gamma(t)}^{r-t}(u^\perp_t),
$$
where $u_t(\tau)=u(t+\tau),\ 0\le\tau\le r-t$. Indeed, the restriction of both $D_vF^r_{q_0}$ and
$D_uF^r_{q_0}$ to the subspace $\{w: w(\tau)=0,\ 0\le\tau\le t\}$ equals $D_{u_t}F^{r-t}_{\gamma(t)}$.
Hence both $D_vF^r_{q_0}(v^\perp)$ and $D_uF^r_{q_0}(u^\perp)$ contain $D_{u_t}F^{r-t}_{\gamma(t)}(u_t^\perp)$.
At the same time, all these three subspaces have codimension 1 in $T_qM$.

Finally, $D_uF^r_{q_0}(u^\perp)={P_t^r}_*D_uF^t_{q_0}(u^\perp)$ and
$D_vF^r_{q_0}(u^\perp)={P_t^r}_*D_{\hat u}F^t_{q_0}(\hat u^\perp).\quad\square$

\begin{lemma}
Let $\varphi:(-1,1)\to M,\ \varphi(0)=\gamma(t)$, be a transversal to $E^t$ smooth curve. Then
any $W^{1,\infty}$-neighborhood of $\gamma|_{[0,t]}$ contains a horizontal curve $\gamma':[0,t]\to M$ such that
$\gamma'(0)=q_0,\ \gamma'(t)\in\varphi((-1,1))\setminus\gamma(t),\ length(\gamma')\le t$.
In particular, $\varphi\bigl((-1,1)\setminus\{0\}\bigr)\cap B_{q_0}(t)\ne\emptyset$.
\end{lemma}

\noindent{\bf Proof.} The statement is local and we can work using particular local coordinates in a neighborhood of $\gamma(t)$. So we may assume that $M=\mathbb R^n=\{(x,y):x\in\mathbb R^{n-1},\ y\in\mathbb R\},\ \varphi(s)=(0,s),\ E^t=\{(x,0):x\in\mathbb R^{n-1}\}$.

Let $V\subset u^\perp$ be an $(n-1)$-dimensional subspace such that $D_uF^t_{q_0}\bigr|_V$ is an isomorphism of $V$ and $E^t_{q_0}$. Let $O_1\subset V$ be the radius 1 Euclidean ball in $V$, where inner product in $V$ is the restriction to $V$ of the $L_2^k[0,t]$ inner product, and $p:(x,y)\mapsto x$ be the projection of $\mathbb R^n$ on $\mathbb R^{n-1}$. We define a map $\Phi^t_\varepsilon:O_1\to\mathbb R^{n-1}$ by the formula:
$$
\Phi^t_\varepsilon(v)=F^t_{q_0}\left(\frac{\sqrt{t}}{\sqrt{t+|\varepsilon v|^2}}(u+\varepsilon v)\right);
$$
then $p\circ\Phi^t_\varepsilon$ is a diffeomorphism of $O_1$ on a neighborhood of 0 in $\mathbb R^{n-1}$ for any sufficiently small $\varepsilon$. We fix such a small $\varepsilon$.

By the Cauchy-- Schwartz inequality, the length of the curve $\tau\mapsto F^\tau_{q_0}(w)$, $0\le \tau\le t,$ doesn't exceed $\sqrt{t}\|w\|,\ \forall w\in\Omega^t_{q_0}$, where $\|\cdot\|$ is the $L_2$ norm. Hence $\Phi^t_\varepsilon(O_1)\in B^t_{q_0}$.

Given $s\in(0,t)$, we consider the map
$$
\Phi^s_\varepsilon:O_1\to\mathbb R^{n-1},\quad
\Phi^s_\varepsilon(v)=F^s_{q_0}\left(\frac{\sqrt{t}}{\sqrt{t+|\varepsilon v|^2}}(u+\varepsilon v)\right).
$$
The maps $\Phi_\varepsilon^s$ uniformly converge to $\Phi^t_\varepsilon$ as $s\to t$. Hence $p\circ\Phi_\varepsilon^s(O_1)$ contains a neighborhood of 0 in $\mathbb R^{n-1}$ for any $s$ sufficiently close to $t$ and $\Phi^s_\varepsilon(O_1)\cap\gamma(-1,1)\ne\emptyset.$
At the same time, $\|(u+ \varepsilon v)|_{[0,s]}\|$ is strictly smaller than
$\|(u+\varepsilon v)|_{[0,t]}\|=\sqrt{t^2+|\varepsilon v|^2}$ if $s<t$, while $\gamma|_{[0,t]}$ is a shortest horizontal path of the length $t$. Hence $\gamma(t)\notin\Phi^s_\varepsilon(O_1).\qquad\square$

We use notations introduced in the proof of Lemma~2. Let
$$
\mathcal B^t=\{w\in\Omega^t_{q_0}:\|w\|^2\le t\},\quad
\mathcal S^t=\{w\in\Omega^t_{q_0}:\|w\|^2=t\}. \eqno (5)
$$
The relation (3) implies that $B_{q_0}(t)=F^t_{q_0}(\mathcal B^t).$

\begin{lemma}
For any sufficiently small connected neighborhood $\mathcal O^t_u$ of $u|_{[0,t]}$ in $\Omega^t_{q_0}$,
there exists $\varepsilon>0$ such that $F^t_{q_0}(\mathcal O^t_u\cap\mathcal B^t)\cap\varphi((-\varepsilon,\varepsilon))$
equals either $\varphi((-\varepsilon,0])$ or $\varphi([0,\varepsilon))$.
\end{lemma}

\noindent{\bf Proof.} Recall that $u$ is a regular point of the map $p\circ F^t_{q_0}\bigr|_{\mathcal S^t}$. According to the standard implicit function theorem, their exists a neighborhood $\mathcal O^t_u$ and a
diffeomorphism $\Psi:\mathcal O^t_u\to\mathbb R^{n-1}\times \mathbb H$, where $\mathbb H$ is a Hilbert space,
such that the following properties hold:
\begin{itemize}
\item $\Psi(u)=(0,0)$.

\item $\Psi(\mathcal O^t_u)=O\times U$, where $O$ is a neighborhood of 0 in $\mathbb R^{n-1}$ and $U$ is a
neighborhood of 0 in $\mathbb H$.

\item $\Psi(\mathcal S^t\cap\mathcal O^t_u)=O\times U_0$, where $U_0=U\cap\mathbb H_0$ and $\mathbb H_0$ is a
codimension 1 subspace of $\mathbb H$.

\item $p\circ F^t_{q_0}\circ \Psi^{-1}(x,w)=x,\quad x\in O,\ w\in U$.
\end{itemize}

Then $\Psi(\mathcal B^t\cap\mathcal O^t_u)=O\times U_-$, where $U_-$ is the intersection of $U$ with a closed half-space in $\mathbb H$ whose boundary is $\mathbb H_0$.

We have:
$$
F^t_{q_0}\circ\Psi^{-1}(x,w)=(x,a(x,w)),\quad x\in O,\ w\in U,
$$
where $a$ is a smooth real function, $a(0,0)=0$. Moreover, Lemma~2 implies that $a(0,U_-)\setminus\{0\}\ne\emptyset$. The subset $a(0,U_-)\subset\mathbb R$ is connected; we have to prove
that 0 is its boundary point.

Assume that it is an interior point: $a(0,w_0)<0,\ a(0,w_1)>0$ for some $w_0,w_1\in U_-$ and consider the broken
line
$$
\bar w=\{sw_0:0\le s\le 1\}\cup\{sw_1:0\le s\le 1\}.
$$
 Then $a(0,\bar w)$ is a segment and 0 is an interior
point of this segment. Hence 0 is an interior point of the segment$a(x,\bar w),\ \forall x$ close to 0. It follows that $\gamma(t)\in int B_{q_0}(t)$, and $\gamma$ is not a shortest horizontal path.$\qquad\square$

We say that $\varphi$ is positively oriented for $\gamma|_{[0,t]}$ if
$$
F^t_{q_0}(\mathcal O^t_u\cap\mathcal B^t)\cap\varphi((-\varepsilon,\varepsilon))=\varphi((-\varepsilon,0])
$$
in the notations of Lemma~3.

\begin{lemma} If $\varphi$ is positively oriented for $\gamma|_{[0,t]}$, then $\varphi$ is positively oriented for any other length-minimizing horizontal path connecting $q_0$ with $\gamma(t)$.
\end{lemma}

\noindent{\bf Proof.} We use notations introduced in Lemma~1. Let $t'>t$ be very closed to $t$ and $\varphi'=p_t^{t'}\varphi$. By continuity, $\varphi'$ is positively oriented for $\gamma|_{[0,t']}$.
Now consider $u|_{[t,t']}$ and apply Lemma~3 to the length-minimizing horizontal path $\gamma|_{[t,t']}$.

Obviously, $F^{t'-t}_{\gamma(t)}\left(\mathcal O^{t'-t}_{u|_{[t,t']}}\cap\mathcal B^{t'-t}\right)\subset
F^{t'}_{q_0}\left(\mathcal O^{t'}_u\cap\mathcal B^{t'}\right)$ Hence $\varphi'$ is positively oriented for
$\gamma|_{[t,t']}$. On the other hand, positive orientation for $\gamma|_{[t,t']}$ is determined by the curve
$\gamma|_{[t,t']}$ and does not depend on the path connecting $q_0$ with $\gamma(t)$.$\qquad\square$

We are now ready to complete the proof of Theorem~1. The set of length-minimizers
$\mathfrak M=\left(F^t_{q_0}\right)^{-1}(\gamma(t))\cap\mathcal B^t$ is strongly compact. Lemmas 3 and 4 imply the existence of a neighborhood $\mathcal O^t_{\mathfrak M}$ and $\varepsilon>0$ such that
$$
F^t_{q_0}(\mathcal O^t_{\mathfrak M}\cap\mathcal B^t)\cap\varphi((-\varepsilon,\varepsilon))=\varphi((-\varepsilon,0]).
$$

It remains to note that $F^t_{q_0}(\mathcal O^t_{\mathfrak M}\cap\mathcal B^t)$ contains a neighborhood of $\gamma(t)$ in $B_{q_0}(t)$. Indeed, it follows from the continuity of $F^t_{q_0}$ in the weak topology, continuity of the subriemannian distance and the fact that weak convergence of a sequence plus convergence of the norms of elements of the sequence imply strong convergence.$\qquad\square$

\medskip
\noindent{\bf Proof of Theorem~2.} Let $\gamma(\tau)=F^\tau_{q_0}(u),\ 0\le\tau\le t,$ be a parameterized
by the length shortest path, $E^t_u=D_uF^t_{q_0}(u^\perp).$

\begin{lemma} If $B_{q_0}(t)$ admits a tangent hyperplane $T_{\gamma(t)}B_{q_0}(t)$, then \linebreak $T_{\gamma(t)}B_{q_0}(t)=E^t_u$.
\end{lemma}

\noindent{\bf Proof.} We use notations (5). For any $\xi\in E^t_u$ there exists a smooth curve and
$s\mapsto w(s)\in\mathcal B_t,\ s\in(-1,1),$ such that $u=w(0),\ \xi=\frac d{ds}F^t_{q_0}(w(s))\bigr|_{s=0}.$
Then $\varphi(s)\doteq F^t_{q_0}(w(s))\in B _{q_0}(t),\ \forall s\in(-1,1),$ and
$\xi=\frac{d\varphi}{ds}(0).$ Hence $\xi\in T_{\gamma(t)}B_{q_0}(t).\qquad \square$

A shortest path $\gamma$ is {\it normal} if $\dot\gamma(t)\notin E^t_u$ and {\it abnormal} if $\dot\gamma(t)\in E^t_u$. Let $\lambda_t\perp E^t_u,\ \lambda_t\in T^*_{\gamma(t)}M\setminus\{0\},$ and $P_{u,\tau}^t:M\to M$ be the diffeomorphism defined by the relation
$P_{u,\tau}^t:q(\tau)\mapsto q(t)$, where $q(\tau),\ 0\le\tau\le t$, is a solution of the differential equation
$\dot q=\sum\limits_{i=1}^ku_i(\tau)X_i(q)$.
We set $\lambda_\tau=P_{u,\tau}^{t*}\lambda_t,\ 0\le\tau\le t,\ \lambda_\tau\in T^*_{\gamma(\tau)}M$.

If $\gamma$ is normal, then $\lambda_\tau,\ 0\le\tau\le t,$ is a solution of the Hamiltonian equation
$\dot\lambda=\vec H(\lambda)$, where $H(\lambda)=\frac 12\sum\limits_{i=1}^k\langle\lambda,X_i(q)\rangle^2,\
\lambda\in T^*_qM,\ q\in M$, and $\vec H$ is the Hamiltonian vector field on $T^*M$ associated to the function $H$.
Moreover, $u_i(\tau)=\langle\lambda_\tau,X_i(\gamma(\tau))\rangle,\ i=1,\ldots,k$. In local coordinates,
$q=x\in\mathbb R^n,\ \lambda=(\xi,x)\in\mathbb R^n\times\mathbb R^n,$
$$
H(\xi,x)=\frac 12\sum\limits_{i=1}^k\langle\xi,X_i(x)\rangle^2, \qquad
\dot\xi=-\frac{\partial H}{\partial x},\ \dot x=\frac{\partial H}{\partial\xi}.
$$
If $\gamma$ is abnormal, then $\lambda_\tau$ satisfies the Goh condition: $\lambda_\tau\perp\Delta^2_{\gamma(\tau)},\ 0\le\tau\le t$.

Lemma 5 and well-posedness of the Cauchy problem for ordinary differential equations imply:

\begin{corollary}
Let $\gamma$ be a normal shortest path. If there exists a one more shortest path (normal or abnormal) connecting $q_0$ with $\gamma(t)$, then $B_{q_0}(t)$ does not admit a tangent hyperplane.
\end{corollary}

Normal shortest paths are projections to $M$ of the solutions of the Hamiltonian system $\dot\lambda=\vec H(\lambda)$; those parameterized by the length are projections of the solutions that belong to the energy level $H(\lambda)=\frac 12$. Moreover, the following fact is well-known (see, for instance, \cite{ABB}):

\begin{prop}
For any compact $K\subset H^{-1}(\frac 12)$ there exists $t'>0$ such that for any solution $\lambda(\tau)\in T^*_{\gamma(\tau)},\ 0\le\tau\le r,$ of the system $\dot\lambda=\vec H(\lambda)$ with the initial condition $\lambda(0)\in K$ the curve $\gamma|_{[0,t']}$ is a unique shortest path connecting $\gamma(0)$ with $\gamma(t')$.
\end{prop}

So short segments of abnormal geodesics are indeed shortest paths.
It may however happen that bigger segments are not like that.

\begin{lemma} Given $t\in(0,r]$, there exists a solution of the system $\dot\lambda=\vec H(\lambda),\ H(\lambda)=\frac 12$, whose projection to $M$ restricted to the segment $[0,t]$ is not a shortest path.
\end{lemma}

\noindent{\bf Proof.} Let $e^{\tau\vec H}:\lambda(\tau)\mapsto\lambda(t)$ be the flow generated by the system
$\dot\lambda=\vec H(\lambda)$ and $\pi:T^*M\to M$ be the standard projection. Assume that
$\tau\mapsto \pi\circ e^{t\vec H}(\lambda_0),\ 0\le\tau\le t$ is a shortest path for any
$\lambda_0\in H^{-1}(\frac 12)\cap T^*_{q_0}M$.

Let $\mu\in\Delta^\perp\setminus(\Delta^2_{q_0})^\perp$; then $H(\mu)=0$ and $H(\lambda_0+\mu)=H(\lambda_0)=\frac 12$.
Consider the family of normal shortest paths
$$
\gamma_s(\tau)=\pi\circ e^{\tau\vec H}(s\mu+\lambda_0),\ 0\le\tau\le t,\ s\in\mathbb R.
$$
Let $\lambda_s(\tau)=e^{\tau\vec H}(\lambda_0+s\mu)$ and 
$$
u_s(\tau)=(\langle\lambda_s(\tau),X_1(\gamma(\tau))\rangle,\cdots,\langle\lambda_s(\tau),X_k(\gamma(\tau))\rangle,
\quad 0\le\tau\le t;
$$
then $\gamma_s(\tau)=F_\tau(u_s),\ \lambda_s(\tau)=P_{u,\tau}^{0*}(s\mu+\lambda_0)$
 and $\lambda_s(\tau)\perp E^\tau_{u_s},\ 0\le\tau\le t$.

Due to the compactness of the set of length minimizers, there exists a sequence $s_n\to\infty\ (n\to\infty)$ such that $u_{s_n}$ strongly converges to a length minimizer $\bar u$ as $n\to\infty$. Then $E^t_{u_{s_n}}$ converges
to $E^t_{\bar u}$ and $P^t_{u_{s_n},\tau}$ converges to $P^t_{\bar u,\tau}$.

We obtain that $\frac 1s\lambda_{s_n}(\tau)=P^{0 *}_{u_{s_n},\tau}\left(\mu+\frac 1{s_n}\lambda_0\right)$ converges    to $\bar\lambda(\tau)=P^{0 *}_{\bar u,\tau}\mu$, where $\bar\lambda(t)\perp E^t_{\bar u}$. The Goh condition implies that $\mu=P^{t*}_{\bar u,0}\bar\lambda(t)$ is orthogonal to $\Delta^2_{q_0}$ that contradicts to our choice of $\mu.\qquad\square$

Let
$$
\gamma(\tau)=\pi\circ e^{\tau\vec H}(\lambda_0),\quad 0\le\tau\le t, \eqno (6)
$$
be a unique horizonal shortest path connecting $q_0$ with $\gamma(t)$. We say that $\gamma$ is {\it extendable} if  there exists $t'>t$ such that $\tau\mapsto\pi\circ e^{\tau\vec H}(\lambda_0),\ 0\le\tau\le t',$ is a unique
horizontal shortest path connecting $q_0$ with $\pi\circ e^{t'\vec H}(\lambda_0)$.

\begin{lemma} The set $\mathfrak L_t$ of all $\lambda_0\in H^{-1}(\frac 12)\cap T^*_{q_0}M$ such that the curve (6) is an extendable shortest path is a nonempty open subset of $H^{-1}(\frac 12)\cap T^*_{q_0}M$.
\end{lemma}

\noindent{\bf Proof.} The set $\mathfrak L_t$ is nonempty because the points connected with $q_0$ by a unique normal shortest path form an everywhere dense subset of $B_{q_0}(r)$ (see\cite{A,ABB}). Let $\lambda_0\in\mathfrak L_t$; we have to prove that a neighborhood of $\lambda_0$ in $H^{-1}(\frac 12)\cap T^*_{q_0}M$ is contained in $\mathfrak L_t$.

First of all, there exists a neighborhood of the curve (6) in $B_{q_0}(r)$ free of endpoints of abnormal shortest paths. This is because the set of these endpoints is closed as the image of a compact set under the endpoint map.

Recall that $\bar t\in(0,t]$ is a {\it conjugate time moment} for the geodesic
$$
\tau\mapsto\pi\circ e^{\tau\vec H}(\lambda_0),\quad 0\le\tau\le t', \eqno (7)
$$
if $\lambda_0$ is a critical point of the map $\lambda\mapsto\pi\circ e^{\bar t\vec H}(\lambda).$
Basic facts (see, for instance, \cite{ABB}):
\begin{itemize}
\item If $(0,t')$ contains a conjugate time moment and $t'$ is not a conjugate time, then the curve (7) is not a shortest path.
\item If any segment of the geodesic has corank 1, then conjugate time moments are isolated.
\end{itemize}
Extendability of the geodesic $\gamma$ implies that $(0,t]$ is free of conjugate time moments.

Further, $\bar t\in(0,t']$ is a {\it cut time moment} for the geodesic (7) if there exists
$\lambda\in H^{-1}(\frac 12)\cap T^*_{q_0}M,\ \lambda\ne\lambda_0$ such that
$\pi\circ e^{\bar t\vec H}(\lambda)=\pi\circ e^{\bar t\vec H}(\lambda_0)$. Then the point
$\pi\circ e^{\bar t\vec H}(\lambda_0)\in B_{q_0}(r)$ is a {\it cut point}. The set of all cut points is the
{\it cut locus}. One more basic fact (see, for instance, \cite{ABB}):
\begin{itemize}
\item Assume that a geodesic
$$
\tau\mapsto\pi\circ e^{\tau\vec H}(\lambda),\quad 0\le\tau\le t', \eqno (8)
$$
does not admit conjugate time moments and endpoints of abnormal shortest paths starting at $q_0$. Then geodesic (8) is a unique shortest path between its endpoint if and only if the interval $(0,t')$ is free of cut time moments.
\end{itemize}

These basic facts and compactness of the set of shortest paths together with Proposition~1 imply that geodesic (8) is a unique shortest path between its endpoints for any $\lambda$ sufficiently close to $\lambda_0$ and $t'$ sufficiently close to $t.\qquad\square$

We are now ready to finish the proof of the Theorem. Lemmas 6, 7 and connectedness of
$H^{-1}(\frac 12)\cap T^*_{q_0}M$ imply that the boundary of $\mathfrak L_t$ is not empty. Let
$\bar\lambda\in\partial\mathfrak L_t$, then
$\tau\mapsto\pi\circ e^{\tau\vec H}(\bar\lambda),\ 0\le\tau\le t,$ is a non-extendable shortest path.

The following known fact completes the proof (see \cite{J,ABB}): the endpoint of a non-extendable shortest
path is either the endpoint of an abnormal shortest path or belongs to the closure of the cut locus.
$\qquad\square$


\begin{thebibliography}{9}

\bibitem{A} A. Agrachev, {\it Any sub-Riemannian metric has points of smoothness.}
Russian Math. Dokl., 2009, v.79, 1--3

\bibitem{ABB} A. Agrachev, B. Barilari, U. Boscain, {\it Introduction to Riemannian and sub-Riemannian geometry.}
http://people.sissa.it/agrachev/books.html

\bibitem{ABC} A. Agrachev, B. Bonnard, M. Chyba, I. Kupka, {\it Sub-Riemannian sphere in Martinet flat case.} J. ESAIM: Control, Optimisation and Calculus of Variations, 1997, v.2,
      377--448

\bibitem{AG} A. Agrachev, J.-P. Gauthier, {\it On the subanalyticity of Carnot-Caratheodory distances.}
Annales de l'institut Henri Poincar\'e--Analyse non lineaire, 2001, v.18, 359--382

\bibitem{CJT} Y. Chitour, F. Jean. E. Tr\'elat, {\it Genericity results for singular trajectories.}
J. Differential Geom., 2006, v.73, 45--73

\bibitem{J} S. Jacquet, {\it Regularity of the sub-Riemannian distance and cut locus.}
Nonlinear control in the year 2000, Vol. 1 (Paris), 521–533, Lecture Notes in Control and Inform. Sci., 258, Springer, London, 2001

\bibitem{T} E. Tr\'elat, {\it Some properties of the value function and its level sets for affine control systems
with quadratic cost.} J. Dynam. Contr. Syst., 2000, v.6, 511--541


\end{thebibliography}
\end{document}